\documentclass[10pt,a4paper]{article}
\usepackage{amssymb}
\usepackage{epsfig}

\newcommand{\Om}{\overline{\Omega}}
\newcommand{\R}{\mathbb{R}}
\renewcommand{\epsilon}{\varepsilon}
\newcommand{\pl}{\Delta_{p}}
\newcommand{\N}{\mathbb{N}}
\newcommand{\cqfd}{\mbox{}\nolinebreak\hfill$\blacksquare$\medbreak\par}
\newtheorem{thm}{Theorem}

\newtheorem{proposition}[thm]{Proposition}

\begin{document}
\title{A note on the moving hyperplane method}\date{}
\author{C. Azizieh and L. Lemaire \footnote{The second author is
supported by an Action de Recherche Concert\'ee de la Communaut\'e Fran\c
caise de Belgique}\\
 D\'epartement de math\'ematique\\ Universit\'e Libre
de Bruxelles\\CP 214 -  Campus Plaine\\ 1050 Bruxelles - Belgique}
\maketitle
Let us consider the problem:
\begin{equation}
\label{problem1}\left\{
\begin{array}{ll}
-\pl u=f(u)&\textrm{in }\Omega,\\
u=0&\textrm{on }\partial\Omega,\\
u\in C^1(\Om), &u>0\textrm{ in }\Omega
\end{array}\right.
\end{equation}
where $1<p\le2$, $\Omega\subset\R^N$ is a bounded convex domain,
$\pl$ is the p-Laplacian operator
defined by $\pl u=\textrm{div}(|\nabla u|^{p-2}\nabla u)$ and
$f:\R\to[0,+\infty)$ is continuous on $\R$, locally Lipschitz
continuous on $(0,+\infty)$ and satisfies
$$
\exists C_0,C_1>0\textrm{ such that } C_0|u|^q\le f(u)\le
C_1|u|^q\quad\forall u\in\R^+
$$
where $q>p-1$.
In \cite{Az}, Ph. Cl\'ement and the first author proved the existence
of a nontrivial positive solution to
(\ref{problem1}) by using continuation methods and establishing
a~priori estimates for the solutions of some nonlinear eigenvalue problem
associated with (\ref{problem1}). The desired a~priori
estimates use a blow up argument as well as some monotonicity and
symmetry results proved by Damascelli and Pacella in \cite{Da} and
generalizing to the p-Laplacian operator with $1<p<2$
the well known results of Gidas--Ni--Nirenberg from
\cite{GNN} and Berestycki--Nirenberg in \cite{BN}.
In their proof, Damascelli and Pacella use a new technique
consisting in moving hyperplanes orthogonal to directions close to
a fixed one. To be efficient, this procedure needs some continuity of some
parameters linked with the moving plane method (see the functions $\lambda_1
(\nu)$ and $a(\nu)$ defined below). Therefore they
assume in their result that $\partial\Omega$ is smooth to insure
this continuity (and only for that reason). However, such a smoothness
hypothesis does not appear
in the case $p=2$ in the classical moving plane
procedure (see \cite{BN}).\\
Our purpose here is to give more precision on the regularity of
the domain $\Omega$  that is needed to have the continuity of the function
$a(\nu)$
and the lower semicontinuity
of $\lambda_1(\nu)$, and so to have the monotonicity and
symmetry
results of \cite{Da}. This question is also important concerning the
existence result from \cite{Az}. Specifically, we ask that the domain be
of class
$C^1$, and we also discuss convexity conditions relating to the continuity of
$\lambda_1(\nu)$.\\

In this paper,
$\Omega$ will
denote an open bounded domain
in $\R^N$ with $C^1$ boundary. We will say that $\Omega$ is
strictly convex if for
all $x,y\in\Om$ and for all $t\in(0,1)$, $(1-t)x+ty\in\Omega$.\\
For any direction $\nu\in\R^N$, $|\nu|=1$, we define
$$
a(\nu):=\inf_{x\in\Omega}x.\nu
$$
and for all $\lambda\ge a(\nu)$,
$$
\Omega_\lambda^\nu:=\{x\in\Omega\,|\,x.\nu<\lambda\},\qquad
T_\lambda^\nu:=\{x\in\Omega\,|\,x.\nu=\lambda\}\,(\ne\emptyset\textrm{ for }
a(\nu)<\lambda<-a(-\nu)).
$$
Let us denote by $R_\lambda^\nu$ the symmetry  with respect to the
hyperplane $T_\lambda^\nu$ and
$$
\begin{array}{rcl}
x_\lambda^\nu&:=&R_\lambda^\nu(x) \,\,\forall x\in\R^N,\\
(\Omega_\lambda^\nu)'&:=&R_\lambda^\nu(\Omega_\lambda^\nu),\\
\Lambda_1(\nu)&:=&\{\mu>a(\nu)\,|\,\forall\lambda\in(a(\nu),\mu),\textrm{
we have (\ref{situation1}) and (\ref{ii})}\},\\
\lambda_1(\nu)&:=&\sup\Lambda_1(\nu)
\end{array}
$$
where (\ref{situation1}), (\ref{ii}) are the following conditions:
\begin{eqnarray} \label{situation1}
{}&(\Omega_\lambda^\nu)' \textrm{ is not internally tangent to }\partial\Omega
\textrm{ at some point }p\notin T_\lambda^\nu&{}\\
\label{ii}
{}&\textrm{for all }x\in\partial\Omega\cap T^\nu_\lambda,\,
\nu(x).\nu\ne0,&{}
\end{eqnarray}
where $\nu(x)$ denotes the inward unit normal to
$\partial\Omega$ at $x$.
Notice that
$\Lambda_1(\nu)\ne\emptyset\quad\textrm{and}\quad\lambda_1(\nu)<\infty$
 since for $\lambda>a(\nu)$
close to $a(\nu)$, (\ref{situation1}) and (\ref{ii}) are satisfied
and $\Omega$ is bounded.
\begin{figure}[h]
\begin{center}
\begin{tabular}{c}
\includegraphics[width=4.7cm,height=3.3cm]{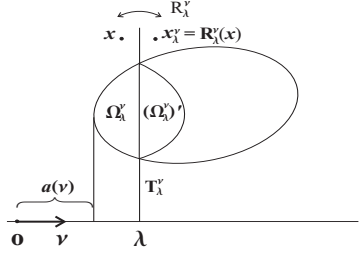}
\end{tabular}
\caption{Illustration of the notations}
\end{center}
\end{figure}

Propositions \ref{luc} and \ref{lem} below give sufficient
conditions on $\Omega$ to guarantee the continuity of the
functions $a(\nu)$ and $\lambda_1(\nu)$, as well as the lower semicontinuity of
$\lambda_1(\nu)$.
\begin{proposition}\label{luc}
Let $\Omega$ be a bounded domain with $C^1$ boundary. Then the function
$a(\nu)$ is continuous with respect to $\nu\in S^{N-1}$.
\end{proposition}
\begin{proposition}\label{lem}
Let $\Omega \subset\R^N$ be a bounded  domain with $C^1$ boundary.
Then the function $\lambda_1(\nu)$ is lower semicontinuous with respect to
$\nu\in
S^{N-1}$. If moreover $\Omega$ is strictly convex, then $\lambda_1(\nu)$
is continuous.
\end{proposition}
As a consequence of these results, we can give more precision on
the conditions to impose to $\Omega$ in the monotonicity result of
\cite{Da}. This result becomes:\vspace*{0.2cm}\\
\textbf{Theorem (Damascelli-Pacella, Theorem 1.1 from
\cite{Da})}\emph{ Let $\Omega$ be a bounded  domain in $\R^N$ with
$C^1$ boundary,
$N\ge2$ and $g:\R\to\R$ be a locally Lipschitz
continuous function. Let $u\in C^1(\bar{\Omega})$ be a weak solution of
$$
\left\{\begin{array}{cl} -\pl u=g(u)&\textrm{in
}\Omega\\u>0&\textrm{in }\Omega,\\ u=0&\textrm{on }\partial\Omega
\end{array}
\right.
$$
where  $1<p<2$. Then, for any direction $\nu\in\R^N$ and for $\lambda$
in the interval $(a(\nu),
\lambda_1(\nu)]$, we have
$u(x)\le u\left(x_\lambda^\nu\right)$ for all $x\in\Omega_\lambda^\nu$.
Moreover $\frac{\partial u}{\partial\nu}(x)>0$ for all
$x\in\Omega_{\lambda_1(\nu)}^\nu\backslash Z$
where $Z=\{x\in\Omega\,|\,\nabla u(x)=0\}$.}\\
\mbox{}\\
Below we prove Propositions \ref{luc} and \ref{lem}
and we give a counterexample of a $C^\infty$ convex
but not strictly convex domain for which $\lambda_1(\nu)$ is not
continuous everywhere.\\
\mbox{}\\
\textbf{Proof of Proposition \ref{luc}:}
Let us fix a direction $\nu\in S^{N-1}$. We shall prove that for all
sequence $\nu_n\to \nu$ with $|\nu_n|=1$, there exists a subsequence still
denoted by $\nu_n$ such that $a(\nu_n)\to a(\nu)$.
Since $\Omega$ is bounded, $(a(\nu_n))$ is also bounded, so passing to an
adequate
subsequence, there exists $\bar a\in\R$ such that
$
a(\nu_n)\to \bar a
$.
We will show that $\bar a= a(\nu)$. Suppose by contradiction that  $\bar
a\ne a(\nu)$. Then either $\bar a<a(\nu)$ or $\bar a>
a(\nu)$.\vspace*{0.2cm}\\
\underline{\textsc{Case 1: $\bar a<a(\nu)$:}} \mbox{ }
Since
$$
a(\nu)=\inf_{x\in\Omega}x.\nu=\min_{x\in\Om}x.\nu=\min_{x\in\partial\Omega}x.\nu
,
$$
there exists $x_n\in\partial\Omega$ such that
\begin{equation}\label{easy}
x_n.\nu_n=a(\nu_n).
\end{equation}
Passing again to a subsequence, there exists $x\in \partial\Omega$ such
that $x_n\to x$ and taking the limit of (\ref{easy}), we get
$
x.\nu=\bar a<a(\nu)
$,
a contradiction with the definition of $a(\nu)$.\vspace*{0.2cm}\\
\underline{\textsc{Case 2: $\bar a> a(\nu)$:}} \mbox{ } There exists
$x\in\partial\Omega$ with $x.\nu=a(\nu)$. For $n$ large,
$|x.\nu_n-x.\nu|=|x.\nu_n-a(\nu)|$ is small, and since $a(\nu_n)\to \bar
a>a(\nu)$, for $n$ large enough we have
$x.\nu_n<a(\nu_n)$, contradicting the definition of $a(\nu_n)$.\cqfd
\mbox{}\\
\textbf{Proof of Proposition \ref{lem}:} We first prove the continuity of
$\lambda_1(\nu)$ if
$\Omega $ is strictly convex. Suppose by contradiction that there exists
$\nu\in
S^{N-1}$ such that $\lambda_1$ is not continuous at $\nu$. Then we can fix
 $\epsilon>0$ and  a sequence $(\nu_n)\subset S^{N-1}$
such that
 $\nu_n\to\nu$ and $|\lambda_1(\nu)-\lambda_1(\nu_n)|>\epsilon$
for all $n\in\N$. Passing to a subsequence still denoted by
$(\nu_n)$, we can suppose that
$$
\textrm{either}\quad\lambda_1(\nu)>\lambda_1(\nu_n)+\epsilon\quad\forall
n\in\N\quad\textrm{or}\quad
\lambda_1(\nu)<\lambda_1(\nu_n)-\epsilon\quad\forall n\in\N.
$$
\underline{\textsc{Case 1:} $\lambda_1(\nu)>\lambda_1(\nu_n)+\epsilon$ for all
$n\in\N$.}\hskip0.6cm
For any fixed $n\in\N$, we have the following alternative:
either there exists $x_n\in  T_{\lambda_1(\nu_n)}^{\nu_n}\cap\partial\Omega $
with $\nu(x_n).\nu_n=0$, or there exists
$x_n\in\left(\partial\Omega\cap\overline
{\Omega_{\lambda_1(\nu_n)}^{\nu_n}}\right)\setminus
T_{\lambda_1(\nu_n)}^{\nu_n}$
with $\left(x_n\right)^{\nu_n}_{\lambda_1(\nu_n)}\in\partial\Omega$.
Passing once again to subsequences, we can suppose that we are in
one of the two situations above for all $n\in\N$. We treat below
each situation and try to reach a contradiction.\vspace*{0.5cm}\\
(1.a) For all $n\in\N$, there exists $x_n\in
T_{\lambda_1(\nu_n)}^{\nu_n}\cap\partial\Omega $
with $\nu(x_n).\nu_n=0$.\vspace*{0.2cm}\\
Passing if necessary to a subsequence, there exist
$\bar\lambda\le\lambda_1(\nu)-\epsilon$
and $x\in T_{\bar\lambda}^\nu\cap\partial\Omega$ such that
$x_n\to x$ and $\nu(x).\nu=0$. This contradicts the definition of
$\lambda_1(\nu)$.\\
\mbox{}\\
(1.b) For all $n\in\N$, there exists $x_n\in\left(\partial\Omega\cap\overline
{\Omega_{\lambda_1(\nu_n)}^{\nu_n}}\right)\setminus
T_{\lambda_1(\nu_n)}^{\nu_n}$
with
$\left(x_n\right)^{\nu_n}_{\lambda_1(\nu_n)}\in\partial\Omega$.\vspace*{0.2cm}\\
Passing if necessary to a subsequence, there exist
$\bar\lambda\le\lambda_1(\nu)-\epsilon$
and $x\in\partial\Omega\cap\overline{\Omega_{\bar\lambda}^\nu}$
such that $x_n\to x$ and $x_{\bar\lambda}^\nu\in\partial\Omega$.
If $x\not\in T_{\bar\lambda}^\nu$, we reach a contradiction with
the definition of $\lambda_1(\nu)$. Suppose now that $x\in
T_{\bar\lambda}^\nu$.
Let us denote $(x_n)^{\nu_n}_{\lambda_1(\nu_n)}$ by $u_n$.
Since $\Omega$ is a $C^1$ domain, it holds that $\nu(u_n).\nu_n\le0$ for all
$n$. By definition of $\lambda_1(\nu_n)$, $\nu(x_n).\nu_n\ge0$. If
$x\in
T_{\bar\lambda}^\nu$, $x=\lim x_n=\lim u_n$ and so $\nu(x).\nu=0$,
which contradicts the definition of $\lambda_1(\nu)$.\\
Observe that we do not use the convexity of the domain in Case
1.\\\mbox{}\\
\underline{\textsc{Case 2:} $\lambda_1(\nu)<\lambda_1(\nu_n)-\epsilon$ for
all $n\in\N$.}:\hskip0.6cm
As in the first case, either there exists $x\in
T_{\lambda_1(\nu)}^\nu\cap\partial\Omega$
with $\nu(x).\nu=0$ or there exists
$x\in\left(\partial\Omega\cap\overline{\Omega_{\lambda_1(\nu)}^\nu}\right)
\setminus T_{\lambda_1(\nu)}^\nu$
such that $x_{\lambda_1(\nu)}^\nu\in\partial\Omega$. We treat the
first situation in (2.a) and the second one in (2.b).\vspace*{0.5cm}\\
(2.a) For $\epsilon$ small enough,
$T_{\lambda_1(\nu)+\frac{\epsilon}{2}}^\nu\cap
\partial\Omega\ne\emptyset$. Since $\Omega$ is strictly convex, there
exists $x'\in
T_{\lambda_1(\nu)+\frac{\epsilon}{2}}^\nu\cap\partial\Omega$ such
that
\begin{equation}\label{eau}
\nu(x').\nu<0.
\end{equation}
For $\epsilon>0$ small enough, there exists $n_0\in\N$ such that for all
$n\ge n_0$,
 the sets $T_{\lambda_1
(\nu)+\frac{\epsilon}{2}}^{\nu_n}\cap\partial\Omega$ are non empty
and since they are compact, we can choose a sequence $(x_n)$
satisfying
$$
   x_n \in T_{\lambda_1
(\nu)+\frac{\epsilon}{2}}^{\nu_n}\cap\partial\Omega, \qquad
|x'-x_n|=\min\left\{|x'-y|\,:\,y\in T_{\lambda_1
(\nu)+\frac{\epsilon}{2}}^{\nu_n}\cap\partial\Omega\right\}.
  $$
Passing if necessary to a subsequence, $x_n\to y$ for some $y\in
T_{\lambda_1
(\nu)+\frac{\epsilon}{2}}^\nu\cap\partial\Omega$
such that
$$
|x'-y|=\lim_{n\to\infty}\textrm{dist}\left( x',T_{\lambda_1
(\nu)+\frac{\epsilon}{2}}^{\nu_n}\cap\partial\Omega\right),
$$
but since this limit is equal to $0$, we infer that
$x'=y$. Now, since $\lambda_1(\nu)<\lambda_1(\nu_n)-\epsilon$ for
all $n\in\N$, $\nu(x_n).\nu_n>0$ for all $n$ and thus
$\nu(x').\nu\ge0$, a contradiction with (\ref{eau}).\\\mbox{}\\
(2.b) The convexity of $\Omega$ implies that
$
x_{\lambda_1(\nu)+\frac{\epsilon}{2}}^\nu\notin\Om.
$
Now,
$
x_{\lambda_1(\nu)+\frac{\epsilon}{2}}^{\nu_n}\to
x_{\lambda_1(\nu)+\frac{\epsilon}{2}}^\nu,
$
so that
\begin{equation}\label{mercure}
x_{\lambda_1(\nu)+\frac{\epsilon}{2}}^{\nu_n}\notin\Om
\end{equation}
for $n$ large enough. But since $x.\nu<\lambda_1(\nu)$ by
definition of $x$, we also have
$
x.\nu_n<\lambda_1(\nu)<\lambda_1(\nu)+\frac{\epsilon}{2}
$
for $n$ sufficiently large, and so
$$x\in\left(\partial\Omega\cap
\Omega_{\lambda_1(\nu)+\frac{\epsilon}{2}}^{\nu_n}\right)\setminus
T_{\lambda_1(\nu)+\frac{\epsilon}{2}}^{\nu_n}$$
for these values of $n$. This fact together with (\ref{mercure})
contradicts the definition of $\lambda_1(\nu_n)$.\\

The proof of the lower semicontinuity follows from Case 1, which
uses only the $C^1$ regularity of the domain.
\cqfd
\mbox{}\\
\textbf{A counterexample in $\R^2$}\\
\mbox{}\\
\begin{figure}[h]
\begin{center}
\begin{tabular}{lc}
{}&{}\\
{}& \vspace{-2cm}$\lambda_1(\nu_n)>\lambda_1(\nu)+\epsilon$\\
\includegraphics[width=5cm,height=2.5cm]{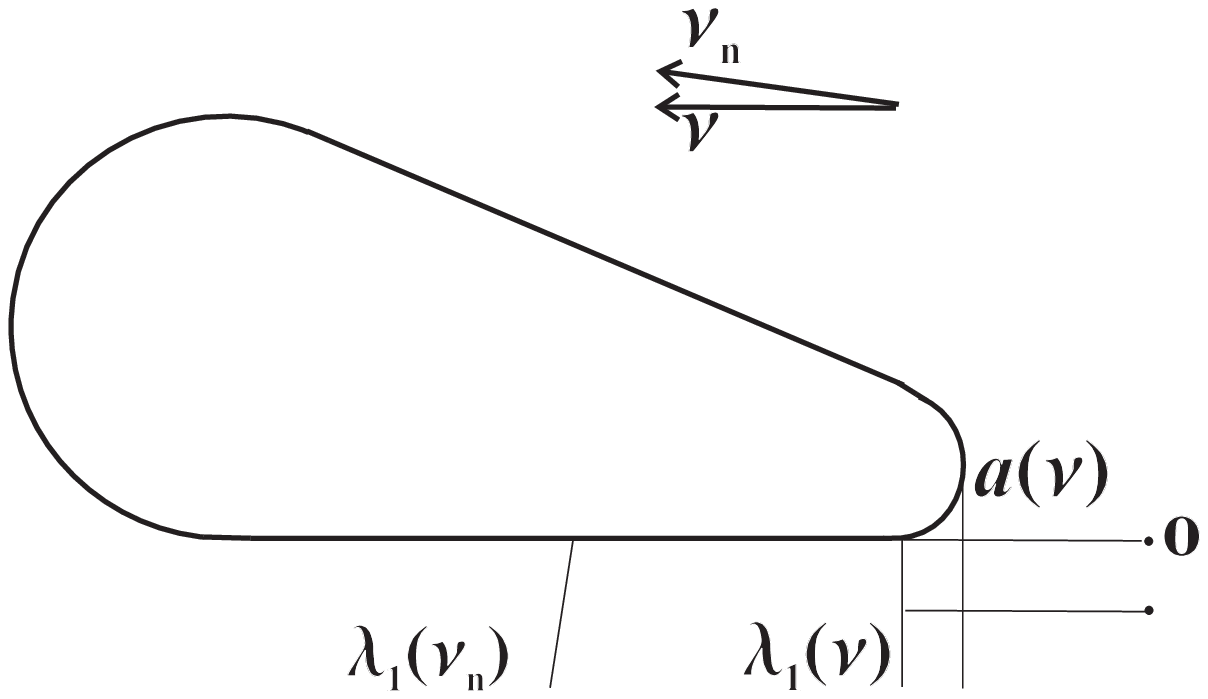}&{ }
\end{tabular}
\end{center}
\caption{Counterexample of a smooth convex but not strictly
convex domain for
which $\lambda_1(\nu)$ is not continuous everywhere.}\end{figure}

This is an example of a convex but not strictly convex domain in
$\R^2$. It contradicts case (2.a) in the proof and indeed, case
(2.a) is the only one using the \emph{strict} convexity. The
example can be made smooth. In fact all is required is a convex
domain in $\R^2$ whose boundary contains a piece of (straight)
line, say  of length $L$. Then for $\nu$ parallel to the line,
there exists a sequence $\nu_n\to\nu$ such that
$\lambda_1(\nu_n)\ge\lambda_1(\nu)
+\frac{L}{2}$.

A variation of this construction will produce similar examples in
higher dimensions.

\end{document}